\documentclass[11pt, twoside]{article}
\usepackage{amssymb}
\usepackage{amsmath, amssymb}\usepackage{cite}\usepackage{mathrsfs}\usepackage[amsmath, thmmarks]{ntheorem}
\usepackage{listings}
\usepackage[titletoc]{appendix}\allowdisplaybreaks
\usepackage{multirow}

\makeatletter
\renewcommand\theequation{\thesection.\@arabic\c@equation}
\makeatother
\newtheorem{thm}{Theorem}[section]%
\newtheorem{lem}[thm]{Lemma}%
\newtheorem{Con}[thm]{Conjecture}%
\newtheorem{Fac}[thm]{Fact}%
\input cyracc.def

\def\f{\noindent}

\def\demo{\f{\bf Proof}\hskip10pt}

\def\qed{\hfill $\Box$}

\begin{document}

\title{{\bf On a conjecture about pattern avoidance of cycle permutations}}
\footnotetext{The work was supported by 2023 Excellent Science and Technology Innovation Team of Jiangsu Province Universities (Real-time Industrial Internet of Things).\\
E-mail addresses: Junyao$_{-}$Pan@126.com}

\author{{\bf Junyao Pan}\\
{\footnotesize Jiangsu Engineering Research Center of Hyperconvergence Application and Security of IoT Devices,}\\
{\footnotesize  Wuxi University, Wuxi, Jiangsu, 214105, P. R. China}
}

\date{}
\maketitle

%\title{\textbf{On the permutations that strongly avoid the pattern $312$ or $231$}}
%\author{Junyao Pan
% \\\\
%School of Sciences, Wuxi University, Wuxi, Jiangsu, 214105, People's Republic of China\\}
%\date {} \maketitle
%
%
%\baselineskip=16pt
%
%
%\vskip0.5cm

\noindent{\small {\bf Abstract:} Let $\pi$ be a cycle permutation that can be expressed as one-line $\pi = \pi_1\pi_2 \cdot\cdot\cdot \pi_n$ and a cycle form $\pi = (c_1,c_2, ..., c_n)$. Archer et al. introduced the notion of pattern avoidance of one-line and all cycle forms for a cycle permutation $\pi$, defined as $\pi_1\pi_2 \cdot\cdot\cdot \pi_n$ and its arbitrary cycle form $c_ic_{i+1}\cdot\cdot\cdot c_nc_1c_2\cdot\cdot\cdot c_{i-1}$ avoid a given pattern. Let $\mathcal{A}^\circ_n(\sigma; \tau)$ denote the set of cyclic permutations in the symmetric group $S_n$ that avoid $\sigma$ in their one-line form and avoid $\tau$ in their all cycle forms. In this note, we prove that $|\mathcal{A}^\circ_n(2431; 1324)|$ is the $(n-1)^{\rm{st}}$ Pell number for any positive integer $n$. Thereby, we give a positive answer to a conjecture of Archer et al.

\vskip0.2cm
\noindent{\small {\bf Keywords:} Pattern avoidance; Cycle permutation; Pell number.


\section {Introduction}

Let $S_n$ denote the symmetric group on $[n]=\{1, 2, \ldots , n\}$. It is well-known that every permutation $\pi$ in $S_n$
can be written either in its cycle form as a product of disjoint cycles or in its one-line notation as $\pi = \pi_1\pi_2\cdots \pi_n$, where $\pi_i = \pi(i)$ for all $i \in [n]$. If $\pi$ is composed of a single $n$-cycle, then $\pi$ is called a \emph{cycle permutation}. Let $\pi= \pi_1\pi_2\cdots \pi_n\in S_n$
and $\tau= \tau_1\tau_2\cdots \tau_k\in S_k$ with $k\leq n$. If there exists a subset of indices $1\leq i_1<i_2<\cdot\cdot\cdot<i_k\leq n$
such that $\pi_{i_s}>\pi_{i_t}$ if and only if $\tau_s>\tau_t$ for all $1\leq s<t\leq k$, then we say that $\tau$ is \emph{contained} in $\pi$
and the subsequence $\pi_{i_1}\pi_{i_2}\cdot\cdot\cdot \pi_{i_k}$ is called an \emph{occurrence} of $\tau$ in $\pi$ and denoted by $\tau\leq \pi$.
For example, $132\leq24153$, because $2,~5,~3$ appear in the same order of size as the letters in $132$. The theory of pattern avoidance in permutations was introduced by Knuth in \cite{K}, which has been widely studied for half a century, refer to \cite{B, V}. A lot of attention has been given to the concept of pattern avoidance over the years. Some interesting and relevant results regarding pattern avoidance can be found in \cite{AE, AG, BC, BS, BD, K, P,P1, SS}.

Let $\pi$ be a cycle permutation in $S_n$. Thereby, $\pi$ can be expressed one-line notation and cycle form as $\pi = \pi_1\pi_2 \cdot\cdot\cdot \pi_n$ and $\pi = (c_1,c_2, ..., c_n)$, respectively. In particular, $\pi$ can also be written $\pi = (c_i,c_{i+1}, ..., c_n,c_1,c_2,...,c_{i-1})$ for each $1\leq i\leq n$; if $c_i=1$ then we call $(c_i,c_{i+1}, ..., c_n,c_1,c_2,...,c_{i-1})$ is the \emph{standard cycle form} of $\pi$. Archer et al. \cite{AB} introduced the notion of pattern avoidance of \emph{one-line and standard cycle form} for a cycle permutation, that is, if $\pi_1\pi_2 \cdot\cdot\cdot \pi_n$ avoids $\sigma$ and $c_ic_{i+1}\cdot\cdot\cdot c_nc_1c_2\cdot\cdot\cdot c_{i-1}$ avoids $\tau$, then $\pi$ avoids $\sigma$ in its one-line form and avoids $\tau$ in its standard cycle form. Archer et al. \cite{AB1} defined the notion of pattern avoidance of \emph{one-line and all cycle forms} for a cycle permutation, namely, if $\pi_1\pi_2 \cdot\cdot\cdot \pi_n$ avoids $\sigma$ and $c_ic_{i+1}\cdot\cdot\cdot c_nc_1c_2\cdot\cdot\cdot c_{i-1}$ avoids $\tau$ for each $1\leq i\leq n$, then $\pi$ avoids $\sigma$ in its one-line form and avoids $\tau$ in its all cycle forms. Let $\mathcal{A}^\circ_n(\sigma; \tau)$ denote the set of cyclic permutations in $S_n$ that avoid $\sigma$ in their one-line form and avoid $\tau$ in their all cycle forms. Archer et al. \cite{AB1} proposed an interesting conjecture about $\mathcal{A}^\circ_n(\sigma; \tau)$, as follows:

\begin{Con}\label{pan1-1}\normalfont([2, Open Questions])
$\Big|\mathcal{A}^\circ_n(2431; 1324)\Big|$ is the $(n-1)^{\rm{st}}$ Pell number.
\end{Con}
In this paper, we prove that the Conjecture\ \ref{pan1-1} is true, and so we obtain the following theorem.
\begin{thm}\label{pan1-2}\normalfont
$\Big|\mathcal{A}^\circ_n(2431; 1324)\Big|$ is the $(n-1)^{\rm{st}}$ Pell number for any positive integer $n$.
\end{thm}

\section {Proof of Theorem 1.2}
It is well-known that the Pell numbers are defined by $P_0=0$ and $P_1=1$, and the recurrence relation $P_n=2P_{(n-1)}+P_{(n-2)}$ for $n\geq2$. Moreover, one easily checks that $\Big|\mathcal{A}^\circ_1(2431; 1324)\Big|=0$ and $\Big|\mathcal{A}^\circ_2(2431; 1324)\Big|=1$ and $\Big|\mathcal{A}^\circ_3(2431; 1324)\Big|=2$ and $\Big|\mathcal{A}^\circ_4(2431; 1324)\Big|=5$. Thus, we see that the Theorem\ \ref{pan1-2} holds for $n=1,2,3,4$. So we shall prove the Theorem\ \ref{pan1-2} by induction on $n$. In other words, it suffices to prove $\Big|\mathcal{A}^\circ_n(2431; 1324)\Big|=2\Big|\mathcal{A}^\circ_{n-1}(2431; 1324)\Big|+\Big|\mathcal{A}^\circ_{n-2}(2431; 1324)\Big|$ for $n\geq5$. Next we start by providing an useful fact that has been pointed out in \cite{AB1}.

\begin{Fac}\label{pan2-0}\normalfont
Let $\pi=(1,c_2, ...,c_{r-1},2,c_{r+1},..., c_n)$ be a cycle permutation in $\mathcal{A}^\circ_n(2431; 1324)$ with $n\geq3$. Then $\{c_2, ...,c_{r-1}\}=\{n-r+3,...,n\}$ and $\{c_{r+1},..., c_n\}=\{3,...,n-r+2\}$. Moreover, if $c_2\neq2$ then the elements after $2$ appear in increasing order.
\end{Fac}
Base on the Fact\ \ref{pan2-0}, we define $$\mathcal{A}^\circ_n(\sigma; \tau)\big|_2^j=\Big\{\pi\in\mathcal{A}^\circ_n(\sigma; \tau)\Big|\pi = (1,c_2, ...,c_{j-1},2,c_{j+1},..., c_n)\Big\}.$$
Thereby, we have
\begin{equation}\label{eq1}
\Big|\mathcal{A}^\circ_n(\sigma; \tau)\Big|=\sum_{j=2}^{n}\Big|\mathcal{A}^\circ_n(\sigma; \tau)\big|_2^j\Big|.
\end{equation}

\begin{lem}\label{pan2-1}\normalfont
Let $n$ be a positive integer with $n\geq4$. Then $\big|\mathcal{A}^\circ_n(2431; 1324)\big|_2^j\big|=\big|\mathcal{A}^\circ_{j-1}(2431; 1324)\big|$ for each $3\leq j\leq n$.
\end{lem}
\demo Consider $\Big|\mathcal{A}^\circ_n(2431; 1324)\big|_2^n\Big|$. For every $(1,c_2, ...,c_{n-1},2)\in\mathcal{A}^\circ_n(2431; 1324)\big|_2^n$, we define a mapping $\mathit{f}$ by the rule that $$\mathit{f}:(1,c_2, ...,c_{n-1},2)\mapsto(1,c_2-1, ...,c_{n-1}-1).$$ We claim that $\mathit{f}$ is a bijection from $\mathcal{A}^\circ_n(2431; 1324)\big|_2^n$ to $\mathcal{A}^\circ_{n-1}(2431; 1324)$. Firstly, we prove that the definition of this mapping is reasonable. Clearly, $(1,c_2-1, ...,c_{n-1}-1)$ avoids $1324$ in its all cycle forms. So it suffices to prove that $(1,c_2-1, ...,c_{n-1}-1)$ avoids $2431$ in its one-line form. Let $(c_1,c_2, ...,c_{n-1},c_n)=\pi_1\pi_2\cdot\cdot\cdot\pi_n$ where $c_1=1$ and $c_n=2$. Note that $\pi_{c_n}=c_{1}$ and $\pi_{c_i}=c_{i+1}$ for $1\leq i\leq n-1$. Pick $\pi'=(\pi_1-1)(\pi_3-1)\cdot\cdot\cdot(\pi_n-1)$. Obviously, $\pi'$ avoids $2431$ in its one-line form. Moreover, we see that $\pi'(1)=\pi_1-1=c_2-1$ and $\pi'(c_i-1)=\pi_{c_i}-1=c_{i+1}-1$ for $2\leq i\leq n-1$. Note $c_{n}-1=1$ and thus $\pi'=(1,c_2-1, ...,c_{n-1}-1)$. Therefore, the definition of this mapping is reasonable. In addition, it is clear that $\mathit{f}$ is an injection, and it can be shown that $\mathit{f}$ is surjection by the same method. Thereby, our claim is true, and so $\Big|\mathcal{A}^\circ_n(2431; 1324)\big|_2^n\Big|=\Big|\mathcal{A}^\circ_{n-1}(2431; 1324)\Big|$.

Consider $\Big|\mathcal{A}^\circ_n(2431; 1324)\big|_2^j\Big|$ for $2<j< n$. By Fact\ \ref{pan2-0}, we see that every $\pi\in\mathcal{A}^\circ_n(2431; 1324)\big|_2^j$ can be expressed as $(1,c_2, ...,c_{j-1},2,3,...,n-j+2)$. For convenience, we set $m=n-j+1$. Now we define a mapping $\mathit{g}$ by the rule that $$\mathit{g}:(1,c_2, ...,c_{j-1},2,3,...,m,m+1)\mapsto(1,c_2-m, ...,c_{j-1}-m).$$ Next we prove that $\mathit{g}$ is a bijection from $\mathcal{A}^\circ_n(2431; 1324)\big|_2^j$ to $\mathcal{A}^\circ_{j-1}(2431; 1324)$. Firstly, we prove that the definition of this mapping is reasonable. Clearly, $(1,c_2-m, ...,c_{j-1}-m)$ avoids $1324$ in its all cycle forms. So it suffices to prove that $(1,c_2-m, ...,c_{j-1}-m)$ avoids $2431$ in its one-line form. Let $(1,c_2, ...,c_{j-1},2,3,...,m,m+1)=\pi_13\cdot\cdot\cdot (m+1)1\pi_{m+2}\cdot\cdot\cdot\pi_n$. Note that $\pi_{c_{j-1}}=2$ and if $c_1=1$ then $\pi_{c_i}=c_{i+1}$ for $1\leq i< j-1$.  Taking $$\pi'=(\pi_1-m)(\pi_{m+2}-m)\cdot\cdot\cdot(\pi_{c_{j-1}-1}-m)(\pi_{c_{j-1}}-1)(\pi_{c_{j-1}+1}-m)\cdot\cdot\cdot(\pi_n-m).$$ Obviously, $\pi'$ avoids $2431$ in its one-line form. Moreover, we see that $\pi'(1)=\pi_1-m=c_2-m$ and $\pi'(c_i-m)=\pi_{c_i}-m=c_{i+1}-m$ for $2\leq i<j-1$. Note $\pi'(c_{j-1}-m)=\pi_{c_j-1}-1=1$ and thus $\pi'=(1,c_2-m, ...,c_{j-1}-m)$. Therefore, the definition of this mapping is reasonable. In addition, it is clear that $\mathit{g}$ is an injection, and it can be shown that $\mathit{g}$ is surjection by the same method. Thereby, $\mathit{g}$ is a bijection and so $\Big|\mathcal{A}^\circ_n(2431; 1324)\big|_2^j\Big|=\Big|\mathcal{A}^\circ_{j-1}(2431; 1324)\Big|$. The proof of this lemma is completed.  \qed

According to equation (2.1) and Lemma\ \ref{pan2-1}, we deduce that for $n\geq5$, $$\Big|\mathcal{A}^\circ_n(2431; 1324)\Big|=\Big|\mathcal{A}^\circ_n(2431; 1324)\big|_2^2\Big|+\sum_{j=2}^{n-1}\Big|\mathcal{A}^\circ_j(2431; 1324)\Big|$$
and $$\Big|\mathcal{A}^\circ_{n-1}(2431; 1324)\Big|=\Big|\mathcal{A}^\circ_{n-1}(2431; 1324)\big|_2^2\Big|+\sum_{j=2}^{n-2}\Big|\mathcal{A}^\circ_j(2431; 1324)\Big|.$$
Thereby, we deduce that $$\Big|\mathcal{A}^\circ_n(2431; 1324)\Big|=2\Big|\mathcal{A}^\circ_{n-1}(2431; 1324)\Big|+\Big|\mathcal{A}^\circ_n(2431; 1324)\big|_2^2\Big|-\Big|\mathcal{A}^\circ_{n-1}(2431; 1324)\big|_2^2\Big|.$$
So far, we have seen that it suffices to prove $$\Big|\mathcal{A}^\circ_n(2431; 1324)\big|_2^2\Big|-\Big|\mathcal{A}^\circ_{n-1}(2431; 1324)\big|_2^2\Big|=\Big|\mathcal{A}^\circ_{n-2}(2431; 1324)\Big|.$$
Inspired by Lemma\ \ref{pan2-1}, we define $$\mathcal{A}^\circ_n(\sigma; \tau)\big|_2^2\big|_3^j=\Big\{\pi\in\mathcal{A}^\circ_n(\sigma; \tau)\big|_2^2\Big|\pi = (1,2,c_3, ...,c_{j-1},3,c_{j+1},..., c_n)\Big\}.$$
Thereby, we have
\begin{equation}\label{eq2}
\Big|\mathcal{A}^\circ_n(\sigma; \tau)\big|_2^2\Big|=\sum_{j=3}^{n}\Big|\mathcal{A}^\circ_n(\sigma; \tau)\big|_2^2\big|_3^j\Big|.
\end{equation}
Next we consider $\mathcal{A}^\circ_n(2431; 1324)\big|_2^2\big|_3^j$ in three situations.

\begin{lem}\label{pan2-2}\normalfont
Let $n$ be a positive integer with $n\geq5$. Then $$\Big|\mathcal{A}^\circ_n(2431; 1324)\big|_2^2\big|_3^n\Big|=\Big|\mathcal{A}^\circ_{n-2}(2431; 1324)\Big|.$$
\end{lem}
\demo For each $(1,2,c_3, ...,c_{n-1},3)\in\mathcal{A}^\circ_n(2431; 1324)\big|_2^2\big|_3^n$, we define a mapping $\mathit{f}$ by the rule that $$\mathit{f}:(1,2,c_3, ...,c_{n-1},3)\mapsto(1,c_3-1, ...,c_{n-1}-1,2).$$ Proceeding as in the proof of Lemma\ \ref{pan2-1}, we deduce that $\mathit{f}$ is a bijection from $\mathcal{A}^\circ_n(2431; 1324)\big|_2^2\big|_3^n$ to $\mathcal{A}^\circ_{n-1}(2431; 1324)\big|_2^{n-1}$. Thereby, $\Big|\mathcal{A}^\circ_n(2431; 1324)\big|_2^2\big|_3^n\Big|=\Big|\mathcal{A}^\circ_{n-1}(2431; 1324)\big|_2^{n-1}\Big|$. It follows from Lemma\ \ref{pan2-1} that $\Big|\mathcal{A}^\circ_n(2431; 1324)\big|_2^2\big|_3^n\Big|=\Big|\mathcal{A}^\circ_{n-2}(2431; 1324)\Big|$, as desired.   \qed

\begin{lem}\label{pan2-3}\normalfont
Let $n$ be a positive integer with $n\geq5$. Then for $3<j<n$, we have $$\Big|\mathcal{A}^\circ_n(2431; 1324)\big|_2^2\big|_3^j\Big|=0.$$
\end{lem}
\demo Suppose $\pi=(1,2,c_3, ...,c_{j-1},3,c_{j+1},...,c_n)\in\mathcal{A}^\circ_n(2431; 1324)\big|_2^2\big|_3^j$ with $3<j<n$. Since $\pi$ avoids $1324$ in its all cycle forms, we infer that $c_{j+1}=4,c_{j+2}=5,...,c_n=n-j+3$. Thereby, $\{c_3, ...,c_{j-1}\}=\{n-j+4,...,n\}$. Let $\pi=\pi_1\pi_2...\pi_n$. Note $\pi_1=2$, $\pi_2=c_3$, $\pi_3=4$ and $\pi_{n-j+3}=1$. Hence, $\pi$ contains $2431$ in its one-line form, a contradiction. Therefore, $\mathcal{A}^\circ_n(2431; 1324)\big|_2^2\big|_3^j=\emptyset$ for $3<j<n$, as desired.   \qed

\begin{lem}\label{pan2-4}\normalfont
Let $n$ be a positive integer with $n\geq5$. Then $$\Big|\mathcal{A}^\circ_n(2431; 1324)\big|_2^2\big|_3^3\Big|=\Big|\mathcal{A}^\circ_{n-1}(2431; 1324)\big|_2^2\Big|.$$
\end{lem}
\demo For every $(1,2,3,c_4, ...,c_{n})\in\mathcal{A}^\circ_n(2431; 1324)\big|_2^2\big|_3^3$, we define a mapping $\mathit{f}$ by the rule that $$\mathit{f}:(1,2,3,c_4, ...,c_{n})\mapsto(1,2,c_4-1, ...,c_{n}-1).$$ Proceeding as in the proof of Lemma\ \ref{pan2-1}, we see that $\mathit{f}$ is a bijection from $\mathcal{A}^\circ_n(2431; 1324)\big|_2^2\big|_3^3$ to $\mathcal{A}^\circ_{n-1}(2431; 1324)\big|_2^{2}$. Thereby, $\Big|\mathcal{A}^\circ_n(2431; 1324)\big|_2^2\big|_3^3\Big|=\Big|\mathcal{A}^\circ_{n-1}(2431; 1324)\big|_2^2\Big|$, as desired.   \qed

According to equation (2.2) and Lemma\ \ref{pan2-2} and Lemma\ \ref{pan2-3} and Lemma\ \ref{pan2-4}, we deduce that $$\Big|\mathcal{A}^\circ_n(2431; 1324)\big|_2^2\Big|-\Big|\mathcal{A}^\circ_{n-1}(2431; 1324)\big|_2^2\Big|=\Big|\mathcal{A}^\circ_{n-2}(2431; 1324)\Big|.$$
Up to now we have completed the proof of Theorem\ \ref{pan1-2}.
%\section{Acknowledgement}
%
%We are very grateful to the anonymous referees for their useful suggestions and comments.

\end{document}